\documentclass{article}
\usepackage{amsmath,latexsym,amssymb,amsthm,enumerate,amscd}
\theoremstyle{plain}
\newtheorem{theorem}{Theorem}[section]

\newtheorem{lemma}[theorem]{Lemma}
\theoremstyle{definition}
\newtheorem{definition}[theorem]{Definition}

\newtheorem{remark}[theorem]{Remark}

\newtheorem*{ack}{Acknowledgment}

\numberwithin{equation}{section}
\numberwithin{table}{section} %TTT
\newcommand{\CI}{CI} 
\newcommand{\Grass}{\mathrm{Grass}}

\newcommand{\PGOR}{\mathrm{PGOR}}
\newcommand{\Levalg}{\mathrm{LevAlg}}
\newcommand{\Hilb}{\mathrm{Hilb}}
\newcommand{\GrAlg}{\mathrm{GrAlg}}
\newcommand{\Points}{\mathrm{LevPoint}}
\newcommand{\Z}{\mathfrak{Z}}
\newcommand{\Ann}{\mathrm{Ann}}
\newcommand{\Soc}{\mathrm{Soc}}
\newcommand{\Gl}{\mathrm{Gl}}
\newcommand{\cha}{\mathrm{char}\ }
\newcommand{\Hom}{\mathrm{Hom}}
\newcommand{\PGor}{\mathrm{PGor}}
\newcommand{\ns}{\footnotesize \it}

\title{Reducible family of height three
level algebras}
\author{Mats Boij\\[.05in]
{\ns Department of Mathematics, KTH, 100 44 Stockholm, Sweden.
}\\[.2in]
Anthony Iarrobino\\[.05in]
{\ns Department of Mathematics, Northeastern University, Boston, MA 02115,
 USA.
}\\[.2in]}

\date{July 14, 2007; revised September 26, 2008}

\begin{document}

\maketitle
%\subjclass{Primary: 13D40; Secondary 14C05}
%\footnote\ddag{AMS 2000 Subj. Class, Primary: 14C05; Secondary: 13D40} 
\begin{abstract}
Let $R=k[x_1,\ldots ,x_r]$ be the polynomial ring in $r$ variables over an infinite field $k$, and let $M$ be the maximal ideal of $R$. Here a \emph{level algebra} will be a graded Artinian quotient $A$ 
of $R$ having socle $Soc(A)=0:M$ in a single degree $j$. The Hilbert function $H(A)=(h_0,h_1,\ldots ,h_j)$ gives the dimension 
$h_i=\dim_k A_i$ of each degree-$i$ graded piece of $A$ for $0\le i\le j$. The embedding dimension of $A$ is $h_1$, and the \emph{type}
 of $A$ is $\dim_k \Soc (A)$, here $h_j$. The family $\Levalg (H)$ of
 level algebra quotients of $R$ having Hilbert function $H$ forms an open subscheme of the family of graded algebras
 or, via Macaulay duality, of a Grassmannian.
\par 
We show that for each of the Hilbert functions $H_1=(1,3,4,4)$ and $H_2=(1,3,6,8,9,3)$ the family $\Levalg (H)$ has several irreducible components (Theorems
\ref{1344thm}\ref{1344thmA}, \ref{36893thmA}).  We  show also that these examples each lift to points. However, in the first example, an irreducible Betti stratum 
for Artinian algebras becomes reducible when lifted to points (Theorem \ref{1344thm}\ref{1344thmB}).
These were the first examples we obtained of multiple components for $\Levalg(H)$ in embedding dimension three. 
\par We show that the second example is the first in an infinite sequence of examples of type three Hilbert functions $H(c)$ 
in which also the number of components gets arbitrarily large (Theorem~\ref{InfiniteTypeThree}).\par
The first case where the phenomenon of multiple components can occur (i.e. the lowest embedding dimension and then 
the lowest type) is that of dimension three and type two.  Examples of this first case have been obtained by the authors (unpublished) and also by J.-O.~Kleppe.  
\end{abstract}

\section{Introduction}
 In Section \ref{sec1.1} we give the context
of this paper and summarize some recent work on level algebras. In
Section
\ref{sec1.2} we discuss the behavior of the graded Betti numbers in a
flat family of Artinian algebras, stating a result of A. Ragusa and A.
Zappal\'{a}, and I. Peeva. In  Section \ref{sec1.3} we describe the
parametrization of
$\Levalg(H)$. In Section \ref{qrsec} we summarize our results.\par
We let $R=k[x_1,\ldots ,x_r]$ be the ring of polynomials in $r$ variables over an infinite field $k$,
 and denote by $M=(x_1,\ldots ,x_r)$ the irrelevant maximal ideal. We will consider graded
standard Artinian algebras
$A=R/I$,
\begin{equation}
A=A_0\oplus A_1\oplus\cdots \oplus A_j, \, A_j\not= 0,
\end{equation}
of socle degree $j$, such that $I\subset M^2$. The \emph{socle} $\Soc
(A)$ is
$\Soc (A)=(0:M)_A=\{ a\in A\mid Ma=0\}$. The Artinian algebra $A$ is 
  \emph{level} if its socle lies in a single degree. Thus, $A$ is level iff
$\Soc (A)=A_j$. We will say that a sequence $H=(h_0,\ldots ,h_j)$ is
\emph{level} if it occurs as the Hilbert function $H(A)$ for some level algebra $A$.  We
denote by $n(A)$, the \emph{length} of $A$, the vector space dimension
$n(A)=\dim_k A$. The \emph{length} $n(\Z)$ of a punctual subscheme $\Z$ of $\mathbb P^n$ is that of a minimal Artinian reduction of its coordinate ring.
The family $\Levalg (H)$ of
 level algebra quotients of $R$ having Hilbert function $H$ forms an open subscheme of the family of graded algebras \cite{Kle1, Kle4, I3}
 or, via Macaulay duality, of a Grassmannian \cite{IK,ChGe,Kle4}.
\subsection{Gorenstein algebras and level algebras.}\label{sec1.1}
When the type is one, the level algebra
is
\emph{Gorenstein}.  Gorenstein algebras occur in many branches of
mathematics, and have been widely studied (see \cite{Mac1,BE,Kle2,Kle3,Hu}).
The level algebras of higher types
$t=2,3,
\ldots $ are a natural generalization of the concept of Artinian Gorenstein
algebras. 
\par To specify a socle degree, type and
Hilbert function of a level algebra A of embedding dimension r is
to specify the Hilbert function and the highest graded Betti number of A,
which must be in a single degree.  Graded Artinian algebras are defined by ideals
$I$ that are intersections of level ideals, so one goal of understanding level
algebras is to be able to extend results about them suitably to
more general Artinian algebras. For example, in \cite{I4} certain results about
level algebras of embedding dimension two are extended to general graded
Artinian algebras of embedding dimension two: in particular the Hilbert functions
compatible with a given type sequence $H(Soc(A))$ are specified. Thus, to understand
level algebras is a step toward the study of the family of Artinian algebras with
given but arbitrary graded Betti numbers.\par For certain pairs
$(r,t)$, there are structure theorems for the level algebras of fixed
embedding dimension $r$ and type $t$. Then the
sequences
$H$ that occur as Hilbert functions for such level algebras
are well understood; and it is known that the family
$\Levalg(H)$ is irreducible, even smooth, of known dimensions. \par In
embedding dimension two, the Hilbert-Burch Theorem states that the defining ideal
of a CM quotient 
$R/I$ is given as the maximal minors of a $(\nu-1)\times \nu$ matrix. When
$H=(1,2,\ldots ,t)$ is a level sequence 
$\Levalg (H)$ is also irreducible and smooth, of known dimension
\cite{I3,I4}.\par
The Buchsbaum-Eisenbud structure
theorem shows that the defining ideal of a height three Gorenstein algebra
$A=R/I$ is generated by the square roots of the
$\nu$ diagonal
$(\nu-1)\times
(\nu-1)$ minors of a
$\nu\times
\nu$ alternating matrix $\mathcal M$, where $ \nu $ is odd: namely $I=$ Pfaffian ideal
of
$\mathcal M$
\cite{BE}. Also in this case, the family
$\PGor(H)$ of Artinian Gorenstein algebras having Hilbert function $H$ is
both irreducible and smooth (see \cite{Di,Kle2}); and the dimension of the
family $\PGor(H)$ is known (several authors: for a survey of these and
related results see
\cite[\S 4.4]{IK}). However, in embedding dimension four there is no structure
theorem for Gorenstein algebras and the family
$\PGor(H)$ is in general neither smooth nor even irreducible
\cite{Bj,IS}.\par
For arbitrary $r$ and types $t$, the maximum and minimum possible level Hilbert functions are well known -- see \cite{I2} for maximum, and \cite{BiGe,ChoI} 
for minimum.
 Recently, several authors have studied
level algebras of small lengths, and types, and embedding dimensions three or
four, with the idea of delimiting the set of possible Hilbert functions that
occur, given the triple
$(h_1,j,t)$ of embedding dimension, socle degree, and type
\cite{GHMS,Za1,Za2}. For example, A. Geramita et al studied type two
level algebras of embedding dimension three, and determined the Hilbert
functions that occur for low socle degree, $j\le 6$. They also gave 
some techniques for determining more generally which $H$ occur in type 2 height
three \cite{GHMS}. Concerning the Hilbert function of height three level
algebras, F. Zanello recently showed that in embedding dimensions at least three,
and for high enough types
$t$, there are level sequences
$H$ that are nonunimodal
\cite{Za2}.  A.~Weiss showed that there are such height three nonunimodal
level $H$ for type $t\ge 5$, and height four nonunimodal level $H$ for $t\ge 3$
\cite{Wei}. Other recent results concern weak Lefschetz properties of level algebras, and liftability of certain level
algebras to points \cite{BjZa,M2,M3,Mmi}.
\par  In embedding
dimension four, M. Brignone and G.~Valla \cite{BriV} noted that there are many
Hilbert function sequences $H=(1,4,\ldots )$, such that the poset
$\beta_{lev}(H)$ of minimal level Betti sequences  that occur for Artinian
algebras
$A\in\Levalg(H)$, has two or more minimal elements. By a result of
A. Ragusa and G. Zappala \cite{RZ} (see Lemma
\ref{RZlem} below), this implies that such families $\Levalg
(H)$ have at least two irreducible components, corresponding to the
minimal elements of
$\beta_{lev}(H)$. \par
 However, there had been no published studies of the component structure
of
$\Levalg(H)$ in embedding dimension three for $t\ge 2$ until we began the present
work. Our first (in the sense of smallest socle
degree, then among those, of smallest length) example of type two height three
$H$ for which $\Levalg(H)$ has two irreducible components is
$H=(1,3,6,10,12,12,6,2)$, which is part of a series we will study elsewhere. Meanwhile,
J.~O. Kleppe has given an example where, proving a
conjecture of the second author about the Hilbert function
$H=(1,3,6,10,14,10,6,2) $, he shows that
$\Levalg (H)$ has at least two irreducible components.  He also notes that by
linking one can construct further such examples \cite[Example 49, Remark
50(b)]{Kle4}.
\subsection{Betti numbers in flat families.}\label{sec1.2}
We here state a preparatory result combining work of several authors.
Recall that the graded Betti numbers for $A$ are defined as follows.
The
minimal resolution of $A$ has the form
\begin{equation}\label{betaeq} 0\rightarrow \bigoplus_k R^{\beta_{r,k}}(-k)
\xrightarrow{\delta_r}\cdots
\xrightarrow{\delta_2}\bigoplus_k R^{\beta_{1,k}}(-k)
\xrightarrow{\delta_1}R\xrightarrow{\delta_0}A\rightarrow 0,
\end{equation}
where the collection $\beta=\left( \beta_{i,j}\right)$ are the graded
 Betti numbers, and the $i$-th (total) Betti number is $\beta_i=\sum_k
\beta_{i,k}$.  The poset $\beta (H)$ has as elements the sequences
$\beta=(\beta_{ik})$ of minimal graded Betti numbers that occur for
Artinian algebras of Hilbert function $H$; the partial order is
\begin{equation}\label{betapoeqn}
\beta\le
\beta' \text{ iff each }\beta_{i,k}\le \beta_{i,k}'. 
\end{equation} 
We denote by
$\beta_{lev}(H)$ the restriction of $\beta (H)$ to the subset of level Betti
sequences
$\beta$: namely those for which 
$\beta_{rk}=0$ except for
$k=j+r$. We denote the corresponding Betti stratum of $\Levalg (H)$ by
$\Levalg_\beta (H)$. A {\it consecutive cancellation} is
when a new sequence
$\beta'$ of graded Betti numbers
 are formed from
$\beta$ by choosing
$\beta'_{ik}=\beta_{ik}$ except that for some one pair $(i_0,k_0)$ of
indices
\begin{equation}
\beta_{i_0,k_0}'=\beta_{i_0,k_0}-1,\text { and }
\beta_{i_0-1,k_0}'=\beta_{i_0-1,k_0}-1.
\end{equation}
The following result is shown by A. Ragusa and G. Zappal\'{a} 
with antecedents by M. Boratynski and S. Greco who showed that the total
Betti numbers are upper semicontinuous on a postulation stratum (\cite{RZ,BG}).
The consecutive cancellation portion was shown by I. Peeva \cite[Remark after
Theorem 1.1]{Pe}, based on a result of K. Pardue \cite{Par}. For further
discussion see
\cite[Theorem 1.1]{M1} and \cite[Remark 7]{Kle4}.
By postulation of a scheme $Z$ we mean the Hilbert function
of its coordinate ring $\mathcal O_Z=R/\mathcal I_Z$.
\begin{lemma}\label{RZlem} Let $Y_T\to T$ be a flat family of punctual
subschemes of
$\, \mathbb P^{r-1}$, such that the postulation $H({\mathcal O}_{Y_t})\mid
t\in T$ is constant, and assume that for $t\in T-t_0$, the
minimal graded Betti numbers of ${\mathcal O}_{Y_t}$ are constant, equal to
$\beta$. Then the minimal Betti numbers $\beta (0)$ at the special point
${\mathcal O}_{Y_{t_0}}$ satisfy $\beta (0)\ge \beta$, in the poset $\beta
(H)$; also $\beta$ may be obtained from $\beta (0)$ by a sequence of consecutive
cancellations.\par
If the poset $\beta_{lev} (H)$ has two incomparable
minimal elements $\beta, \beta'$, then $\Levalg (H)$ has at least two
irreducible components corresponding to closures of open subfamilies of
$\Levalg_\beta (H)$, $
\Levalg_{\beta'}(H)$.
\end{lemma}
\subsection{Parametrization of the family $\Levalg (H)$}\label{sec1.3}
We let $H= (1,h_1,h_2,\ldots h_j)$, with $h_1=r$ be a fixed level sequence. For
a vector subspace $V\subset R_i$ we denote by $R_aV$ the vector subspace
\begin{equation*}
R_aV=\langle fv, f\in R_a, v\in V\rangle \subset R_{a+i}.
\end{equation*}
We denote by $\mathcal R$ the ring of divided powers for $r$ variables over $k$, upon which
$R$ acts by differentiation or ``contraction'' (see \cite[Appendix A]{IK}, this is an avatar of
Macaulay duality \cite{Mac1}). We will
use the ring
$\mathcal R$ in the second parametrization $L(H)$ for level algebras in B. just below. We
denote by
$r_i=\dim_kR_i= \binom{i+r-1}{i}$.\par There are two natural ways to parametrize
the family of level algebra quotients of $R$ having Hilbert function $H$:
\begin{enumerate}[A.]
\item\label{schA}  $\Levalg(H)$ is an open subscheme of $\GrAlg (H)$, the family
of graded algebras quotients of $R$ having Hilbert function $H$. The family
$\GrAlg (H)$ is the closed subscheme of 
\begin{equation}
\prod_{i=2}^j \Grass (r_i-h_i,r_i), 
\end{equation}
parametrizing those sequences of subspaces
\begin{align}\label{parameq}
&(V_2,\ldots ,V_j), V_i\subset R_i, \dim V_i=r_i-h_i, \text { such that }\notag\\
&R_1V_2\subset V_3,
\ldots ,R_1V_{j-1}\subset V_j.
\end{align}
The condition \eqref{parameq} is just that the sequence of vector spaces
forms the degree two to $j$ graded components of a graded ideal of $R$. 
For further detail see \cite{Kle1,Kle3,Kle4}, or the discussion of
the ``postulation Hilbert scheme'' in
\cite{IK}, or \cite[Definition 1.9]{I3}.
\item\label{schB} Closed points $p_A, A=R/I$ of $\Levalg(H)$ correspond by Macaulay duality 1-1
to vector subspaces $\mathcal W=(I_j)^\perp$ of
$\mathcal R_j=\Hom(R_j,k)$, hence to points of the Grassmannian $\Grass(h_j,\mathcal R_j)$.
We define the locally closed subscheme $L(H)\subset \Grass(h_j,\mathcal R_j)$ by
the ``catalecticant'' conditions specifying 
\begin{equation*}
\dim R_{j-1}\circ \mathcal W=h_1, \ldots
,\dim R_{1}\circ \mathcal W=h_{j-1}.
\end{equation*}
This parametrization is introduced in the type one Gorenstein case in
\cite[Section 1.1]{IK}, and for general level algebras in
\cite{ChGe}.
\end{enumerate}
The closed points of the two parametrizations $\Levalg(H)$ and $L(H)$ are evidently
the same (in
$\cha k=p$ $\le j$ one must use the ``contraction action'' of $R_i$ on $\mathcal R_j$). By the
universality property of the family of graded algebras \cite{Kle1} there is a
morphism $\iota: L(H)_{red}\to \Levalg(H)$, where $L(H)_{red}$ denotes the
reduced scheme structure (see \cite[Problem 12]{Kle3}). Recently J.-O. Kleppe has
shown that $L(H)$ and $ \Levalg(H)$ give the same topological structures
\cite[Theorem 44]{Kle4}, extending his earlier result that there is an
isomorphism between the tangent spaces to $\Levalg (H)$ and to $L(H)$ at
corresponding closed points.
\subsection{Questions and Results}\label{qrsec}\par
What is a good description of $\Levalg(H)$? From the point of view of deformations we should
answer as fully as possible:\par
\begin{enumerate}[i.]
\item What are the possible level Betti sequences $\beta$ compatible with $H$?
\item What is the dimension of each Betti stratum $ \Levalg_\beta (H)$? 
\item What is the closure of $\Levalg_\beta (H)$?
\item\label{qiv} What are the irreducible components of $\Levalg (H)$ and of $\Levalg_\beta (H)$?
\begin{enumerate}[a.]
\item\label{qva} Can the component structure of $\Levalg(H)$ be related to that of
an appropriate Hilbert scheme of points or of curves on $\mathbb P^{r-1}$?
\item\label{qvb} Do the components of $\Levalg (H)$ lift to components of $\Points
(T)\subset \mathbb P^r,
\Delta T=H$, and further to families of curves?
\end{enumerate}
\end{enumerate}
Of course, in the absence of a structure theorem (so when $(r,t)\not= (2,t)$ or $(3,1)$), it is hopeless to answer all 
of these questions for all $H$, even in
embedding dimension three.  However, we can answer them for certain $H$, that might either be of special interest or
suggest patterns that are frequent.  A productive approach has been that of
Question (\ref{qva}.)
\cite{Bj,IK,Kle3,Kle4}). J.-O.~Kleppe in \cite{Kle4} establishes in many cases
a 1-1 correspondence between the set of irreducible components of $\Levalg (H)$
and those of a suitable Hilbert scheme of points.
\par\smallskip
\subsubsection{Results.}
In this article we first answer 
the questions above for $H_1=(1,3,4,4)$ and
$H_2=(1,3,6,8,9,3)$, perhaps the simplest cases in height three where $\Levalg(H)$
has several components. \par
We show that for each Hilbert function $H=H_1$ or $H=H_2$
the family
$\Levalg (H)$ has several irreducible components (Theorems
\ref{1344thm}\ref{1344thmA},
\ref{36893thmA}).  We determine the Betti strata and their closures for
each of these two Artinian examples. We then show that each of these examples lifts
to families of points having several components (Theorems
\ref{1344thm}\ref{1344thmB},
\ref{36893thmB}). However, in the first example, an irreducible Betti stratum for
Artinian algebras becomes reducible when lifted to points, although the
overall component structures for $\Levalg (H_1)$ and $\Points (T_1)$ correspond.
\par These two families $\Levalg(H)$ have different
behavior with respect to the poset
${\beta}_{lev}(H)$ of graded level Betti number sequences compatible with $H$. In
the first, $H_1 = (1,3,4,4)$, there are two minimal elements of $\beta (H_1)$  and
the Ragusa-Zappal\'{a} result applies. In the other, $H_2=(1,3,6,8,9,3)$ there is
a unique minimum element of
$\beta_{lev}(H)$, and we must use a different argument to show the
reducibility of $\Levalg (H)$. \par 
We show in Section \ref{infsec} that the example $H_2$ is the first in an infinite series 
of height three Hilbert functions of type three, $H(c)$, $c\ge 3$, where $H(c)$ has socle degree $2c-1$ and satisfies
$$
H(c)_i = \min \{r_i-2r_{i-c}, 3r_{2c-1-i}\}, \qquad 0\le i\le 2c-1;
$$
and such that the number of irreducible components of $\Levalg ( H(c))$ is bounded from below by $(1-1/\sqrt2)c$ (Theorem~\ref{InfiniteTypeThree} and
 Remark \ref{endrmk}).
This result uses the connection mentioned above in Question (iva.) to the Hilbert scheme of points on $\mathbb P^2$.
We denote by $G(c)=\Grass(2,R_c)$ the Grassmanian parametrizing pencils of degree $c$ plane curves; we let $X(c)$ denote the closed subset of $G(c)$ that is the complement of the
open dense set of pencils $\langle f,g\rangle$ spanned by a CI. Then $X(c)$ is the union of irreducible components $G(c)_a, 1\le a \le c-1$ corresponding to the degree
$a$ of the base component of the pencil.
A
delicate issue is whether the coordinate ring of the variety that is a union of the base curve $C_a$ and a general enough
 complete intersection of bidegree $(c-a,c-a)$, has type three Artinian quotients of
Hilbert function $H(c)$:  we show this using the uniform position property of general enough complete intersections (Lemma \ref{lemmaHF}).
\par
We plan to study elsewhere further series of height three 
Hilbert functions, for which $\Levalg(H)$ has several irreducible components; one such series begins with 
 $H=(1,3,6,10,12,12,6,2)$ of type two, but at the time of writing we have been not able to show that
the series is infinite.
\section{Families $\Levalg (H)$ having several irreducible components}
We now state and prove our main results, outlined above.  Henceforth we let $R=k[x,y,z]$ so $r_i= \binom{i+2}{i}$. Note that $k$ is algebraically closed in Section \ref{344},
infinite in Section \ref{6893}, and algebraically closed of characteristic zero in Section \ref{infsec}.

\subsection{The family $\Levalg(H_1), H_1=(1,3,4,4)$.}\label{344}
 We first consider the Artinian case. Let $R=k[x,y,z]$ with $ k $ algebraically closed, and
consider
$
\Levalg(H_1), H_1=(1,3,4,4)$.  The subfamily $C_1$ parametrizes type
4, socle degree three level quotients $A$ of complete intersections $B=R/J$
where $J=(f_1,f_2)$ has generator degrees (2,2). Thus, $A=R/I$, where the
defining ideal satisfies (here $m=(x,y,z)$)
\begin{equation}
I=(f_1,f_2,m^4). 
\end{equation}
As a variety, $C_1\cong\CI (2,2)$, the variety parametrizing complete
intersections $B=R/J$, of generator degrees $(2,2)$ for $J$; evidently
$C_1$ may be regarded as a dense open subset of $\Grass(2,R_2)$, so satisfies
\begin{equation*}
\dim C_1=2\cdot 4=8.
 \end{equation*}
It is easy to check that the monomial ideal
\begin{equation}\label{ciex}
I(1)=(x^2,y^2, m^4)
\end{equation}
defines $A(1)=R/(I(1))$ in this component, and has the minimal
resolution $\beta (1)$ of all general enough elements of $C_1$, given in Table
\ref{CI1table} (we use the standard ``{\sc{Macaulay}}'' notation). Note that it is not possible to split off the redundant
term $R(-4)$ in the minimal free resolution $\beta(1)$ of $A(1)$, as the two quadrics need such a syzygy: thus $\beta (1)$ is a minimal
element of the poset $\beta (H)$. 
\begin{table}[htb]
\begin{center}
\begin{tabular}{|l| r r r r |}
\hline
total&       1& 6 & 9  & 4 \\ \hline 
    0: &     1&   - & -  &  - \\
    1: &    - &   2& -&    - \\
    2: &    - &   -  &  1 &   -\\
    3:&     - &   4&  8&    4  \\\hline

\end{tabular}
\caption{Graded Betti numbers $\beta (1)$ for $H_1=(1,3,4,4)$, CI related
Artinian algebra.}\label{CI1table}
\end{center}
\end{table}
\par\smallskip
The subfamily $C_2$ parametrizes type 4 socle degree three
quotients $A=R/I$ of the coordinate ring $B$ of a line union a point
in
$\mathbb P^2$: that is, we let $B=R/(I_2)$, where $
I_2\cong
\langle xy,xz\rangle$ or $I_2\cong \langle x^2, xy \rangle$ up to
$\Gl(3)$ linear map (coordinate change). Then
$R_1I_2$ has vector space dimension five, and also codimension five in $R_3$,
because of the one linear relation on $I_2$. We have
\begin{equation}
\dim C_2= 8. 
\end{equation} This is
$2+2$ for the choice of a point and a line, plus $1\cdot 4$ for the
choice of a cubic form, an element of the quotient vector space $R_3/R_1I_2$. 
Thus, for $A=R/I\in C_2$ we have, after a linear coordinate change,
\begin{equation}
I\cong(xy,xz,f, W), f\in R_3, W\subset R_4, \dim_k W=3,
\end{equation}
or similar generators with $I_2=(x^2,xy)$.
A monomial ideal $I(2)$ determining an Artinian algebra $A(2)$ in $ C _2$ is
\begin{equation}\label{lineeq}
I(2)=(x^2,xy, z^3, y^4, y^2z^2,y^3z),
\end{equation} whose minimal
resolution $\beta (2)$ is that of Table \ref{Cubictable}. Again, although there is a redundant term $R(-3)$ in the minimal free resolution, it cannot be
split off, so $\beta (2)$ is also a minimal element of the poset $\beta (H)$.
\begin{table}[htb]
\begin{center}
\begin{tabular}{|l| r r r r |}
\hline
total&       1& 6 & 9  & 4 \\ \hline 
    0: &     1&   - & -  &  - \\
    1: &    - &   2& 1&    - \\
    2: &    - &   1  &  - &   -\\
    3:&     - &   3&  8&    4  \\\hline

\end{tabular}
\caption{Graded Betti numbers $\beta (2)$ for $H_1=(1,3,4,4)$,
Artinian algebra related to a line.}\label{Cubictable}
\end{center}
\end{table} 
\begin{table}[htb]
\begin{center}
\begin{tabular}{|l| r r r r |}
\hline
total&       1& 7 & 10  & 4 \\ \hline 
    0: &     1&   - & -  &  - \\
    1: &    - &   2& 1&    - \\
    2: &    - &   1  &  1 &   -\\
    3:&     - &   4&  8&    4  \\\hline

\end{tabular}
\caption{Graded Betti numbers $\beta (3)$ for $H_1=(1,3,4,4)$,
Artinian algebras in $\overline{C_1}\cap \overline{C_2}$.}\label{Intersecttable}
\end{center}
\end{table}\par
We now show that there is only one other Betti sequence possible for Artinian 
algebras in $\Levalg(H_1)$, namely $\beta (3)$, the supremum in $\beta(H_1)$ of
$\beta (1)$ and $\beta (2)$ (see Table \ref{Intersecttable}). A monomial ideal
$I(3)$ defining an algebra $A(3)=R/I(3)$ in $C_3:= \Levalg_{\beta (3)}(H_1)$  is
\begin{equation}\label{eqintideal2}
(x^2,xy,y^3,xz^3,y^2z^2,yz^3,z^4). 
\end{equation}
Since we show
in Theorem \ref{1344thm} that the Betti stratum $C_3=\overline{C_1}\cap\overline{C_2}$, the algebra $A(3)$ is a
simple example of an obstructed level algebra. \par\noindent
\begin{lemma}\label{bettistratumlemma}
The Betti stratum $C_3=\Levalg_{\beta
(3)}(H_1)$ is irreducible of dimension seven. There are no further Betti
strata in $\Levalg (H_1)$, other than $\beta (1),\beta (2),$ and $\beta
(3)$.
\end{lemma}
\begin{proof} Let $A=R/I$ be a level algebra of Hilbert function $H(A)=H_1$. In
order to have a linear relation among the two degree two generators of $I$, we
must have
\begin{equation}\label{i2eq}
I_2\cong
\langle xy ,xz \rangle,\text{ or }I_2\cong \langle x^2,xy\rangle,
\end{equation}
 up to
isomorphism. \par
We now prove that $\dim \Levalg_{\beta
(3)}(H_1)=7$.
Suppose that we are in the first case of \eqref{i2eq}. In
order to have a further relation in degree four among the degree two
generators $xy,xz$
and the degree three generator
$f$, we must have
\begin{equation}\label{relationeq}
\ell \cdot f= q_1\cdot xy+q_2\cdot xz, \text { where }\ell\in R_1, \,\,
q_1,q_2\in R_2,
\end{equation}
whence $x$ divides $\ell$ or $f$. If $x$
divides $f$, then $f$ mod $(xy,xz)$ satisfies
$f=ax^3, a\in k^\ast$, and it follows that
$R_1\cdot x^2\subset I$, implying $x^2\in \Soc(A)$, contradicting the assumption
that
$A$ is level of socle degree three. Thus,
$ x$ divides $\ell$, and $f\in (y,z)$. Since we may now assume that $f$,
which may be taken mod $(xy,xz)$, has no terms involving $x$, we have that
$f\in R'_3$ where $R'=k[y,z]$. Since $k$ is algebraically closed, up to
an isomorphism of
$R'$, we may assume one of
$f=y^3$, $f=y^2z$, or $f=yz(y+z)$, the last being the generic case,
specializing to the two others. The dimension count for this
subfamily is seven: two for the choice of $x\in R_1$, two for the choice of
$\langle y,z\rangle\subset R$, determining $R'$, and three for the choice
of an element $f$ of $R'_3=\langle y^3,y^2z,yz^2,z^3\rangle $ mod
$k^\ast$. \par Next let us assume that $I_2=\langle x^2,xy\rangle$, the second
case of \eqref{i2eq}. Then \eqref{relationeq} is replaced by
\begin{equation}\label{specialeq}
\ell \cdot f= q_1\cdot x^2+q_2\cdot xy, \text { where }\ell\in R_1, \,\,
q_1,q_2\in R_2.
\end{equation}
If $x$ divides $f$, then $f$ mod
$(x^2,xy)$ satisfies $f=axz^2, a\in k^\ast$, whence
$R_1\cdot xz \subset I$ contradicting that
$A
$ is level. Hence $x$ divides $\ell$,  $f\in (x,y)$, and $f$ mod
$(x^2,xy)$ satisfies
\begin{equation}
f\in \langle
xz^2,y^3,y^2z,yz^2\rangle.
\end{equation}
 The dimension count for this
subfamily (where $I_2\cong \langle x^2,xy\rangle$), is six: two for the
choice of  $x\in R_1$, one for the choice of $\langle x,y\rangle$
containing $x$, and three for the choice of $f$ up to scalar from a 
four dimensional vector space.  This completes the proof that  $\dim \Levalg_{\beta
(3)}(H_1)=7$.
\par
We now show that $C_3=\Levalg_{\beta (3)}(H_1)$ is irreducible.
Assume that $A=R/I\in \Levalg_{\beta
(3)}(H)$, and that $\ell=x$ in \eqref{specialeq}. We deform
$(x^2,xy)$ to $(x(x+t(z-x),xy), t\in k$, and correspondingly deform $f$ to $f(t)$
satisfying
\begin{equation*}
f(t)= q_1\cdot (x+t(x-z))+q_2\cdot y,
\end{equation*}
 Then the relation
\eqref{specialeq} deforms
to
\begin{equation*}
x \cdot f(t)= q_1\cdot (x^2+tx(x-z))+q_2\cdot xy.
\end{equation*}
This gives a deformation of $A=R/I$ to $A(t)=R/I(t)$ whose fiber over $t\not=0$
is an algebra with $I(t))_2\cong \langle xy,xz\rangle$. Defining for $t\in k$,
$I(t)=(x(x+t(x-z)), xy,f(t),m^4)$, we have constant length, hence flat, family
of Artinian algebras. Thus $C_3$ is irreducible.
\par
We next show that there are no further level Betti sequences for $H_1$. Suppose by
way of contradiction that $A=R/I\in \Levalg (H_1)$ and that there are two or more
relations in degree four, among the three generators of the ideal $I$ having
degrees two and three. More than two is not possible by Macaulay's theorem. If
there are exactly two such relations, the algebra
$A'=R/(I_{\le 3})$ would have Hilbert function  $H(A')=(1,3,4,4,5,\ldots
)$; the growth of $H(A')$ from 4 to 5 in degrees three and four is the maximum
possible by Macaulay's Hilbert function theorem \cite{Mac2},(see \cite[p.
155-156]{BH}) since
$4_3=\binom{4}{3}$ so the maximum possible for $H(A')_4$ is $4_3'=\binom{5}{4}$.
It follows from the Gotzmann Hilbert scheme theorem
(\cite{Gotz}, see also \cite{IKl}) that
$I_3=(I_\Z)_3$, where $\Z$ is a projective variety with Hilbert polynomial
t+1, and regularity degree one, the number of terms in the Macaulay
expansion. Thus,
$\Z$ is a line. It follows from $I$ level that
$I_{\le 3}=(I_\Z)_{\le 3}$, implying $H(A)=(1,2,3,4)$, a contradiction. 
We have shown that there is at most one relation in degree four among the
generators of $(I_{\le 3})$.\par
 Since
$A=R/I\in \Levalg (H_1)$ implies that the minimal generators of $I$ have
degrees at most four, we have ruled out any Betti sequence greater than
$\beta (3)$ for level algebras of Hilbert function $H$. This completes
the proof of Lemma \ref{bettistratumlemma}.
\end{proof}
\par
\begin{remark}
It is easy to see that here as for most $H$, $\beta (H)\not=\beta_{lev}(H)$,
as the known maximal graded Betti numbers are rarely level. Here, taking $I=\Ann
(X^3,Y^3,(X+Y)^3,(X+2Y)^3,Z^2)$, with $X,Y,Z\in S$ gives such a maximum $\beta$ in $\beta (H_1)-\beta_{lev}(H_1)$.
\end{remark}
The first component of $\Levalg(H_1)$ in the Theorem below
is $\overline{C_1}$, the closure in $\Levalg(H_1)$ of the subfamily $C_1$ for which
$I_2$ defines a complete intersection. Recall that the algebras in $C_1$
comprise the Betti stratum
$\beta (1)$ of Table
\ref{CI1table}. The second component is
$\overline{C_2}$, the closure in $\Levalg (H_1)$ of the subfamily $C_2$ for
which
$I_2$ determines a point union a line. Recall that the algebras in $C_2$ comprise
the Betti stratum
$\beta (2)$ of Table
\ref{Cubictable}, and that those in $C_3$ comprise the $\beta (3)$ stratum of $\Levalg(H_1)$. \par
 We let $T_1=(1,4,8,12,12,\ldots )$ with first difference $\Delta T_1=H_1$, and now define corresponding
subfamilies
$C_1',C_2', C_3'$ of
$\Points (T_1)$, the family of length $n(\Z)=12$ punctual subschemes $\Z$ of $\mathbb P^3$ having
postulation $T_1$ (equivalently, having $h$-vector $H_1$).  We let $C_1'$
parametrize the punctual schemes $\Z$ lying over $C_1$, in the sense that
some minimal reduction of
$R/I_\Z $ (some quotient $A$ of $R/I_\Z$ by a linear element
$\ell\in A_1$, such that  $\dim_k A=n(\Z)$) is an Artinian algebra in $C_1$.
We similarly define $C_2'$ lying over $C_2$, and $C_3'$ over $C_3$. \par We
denote by $D_a'$,$D_b'$ and $D_{ab}'$ the following three subfamilies of $C_3'$.
In describing, say, eight points on a plane cubic, we allow certain degenerate
configurations -- as a length 8 punctual scheme on the cubic --- provided that
the resulting punctual scheme remains in $C'_3$: that is, it has an Artinian quotient that is in $C_3$. In particular the Artinian quotient must have the Hilbert function $H_1$, a condition which limits the
amount of degeneracy. (We don't here attempt to delimit the allowable degenerations). We let
$D_a'$ parametrize the subfamily of punctual subschemes of $\mathbb P^3$ that are unions of nine
points on an irreducible plane
cubic, and three points on a line meeting the cubic. Let 
$D_b'$ be the subfamily parametrizing subschemes that are unions of eight general enough
points on a plane and four points on a line not in the plane. We let $D_{ab}'$
parametrize general enough subschemes that are unions of nine points on a plane cubic, of which
one is the intersection of the line and cubic, union three more points on the
line; or a degeneration with a double point at the intersection of the line and
cubic.
\begin{theorem}\label{1344thm} \enumerate[A.]
\item\label{1344thmA} $\Levalg(H_1), H_1= (1,3,4,4)$ has the
two irreducible components $\overline{C_1},\overline{C_2}$  defined above, each
the closure of a Betti stratum.  Both components have dimension eight. We have
\begin{equation}
\overline{C_1}\cap \overline{C_2}=C_3:=\Levalg_{\beta
(3)}(H_1).
\end{equation} 
\item\label{1344thmB}
$\Points (T_1), T_1=(1,4,8,12,12,\dots )$ has two irreducible components, each the
closure of a subfamily
$C_1',C_2'$ lying above $C_1,C_2$, respectively. Algebras in $C_1'$ are the
coordinate rings of schemes comprising 12 points lying on the complete
intersection curve of two quadric surfaces in
$\mathbb P^3$. Elements of $C_2'$ parametrize (an open subfamily of)
schemes comprised of nine points lying on a plane cubic, union  three points on a
line in $\mathbb P^3$. Each subfamily $C_1'$ and $C_2'$ has dimension 28. The
Betti stratum $C_3'$ lies in $\overline{C_1'}$ and has two irreducible
components, the closures of
$D_a'$ and
$D_b'$, each of dimension $27$. The intersection
\begin{equation}
\overline{C_1'}\cap \overline{C_2'}=\overline{D'_a}\subset C_3'.
\end{equation}  
\par The intersection of $\overline{D_a'}$ and $\overline{D_b'}$ is $D_{ab}'$, and
has dimension $26$.
\end{theorem}
\begin{proof} To show that $\Levalg(H_1)$ has two irreducible components, it
suffices by the Ragusa-Zappal\'{a} Lemma \ref{RZlem} to note that the Betti
resolutions of Tables
\ref{CI1table} and
\ref{Cubictable}, which occur for the monomial ideals given in
\eqref{ciex},\eqref{lineeq}, are each minimal in $\beta (H_1)$. This is true as the
two degree
 two generators of $I$ must have a
quadratic relation as in Table \ref{CI1table} unless they have a linear
relation as in Table~\ref{Cubictable} (or see Lemma
\ref{bettistratumlemma}.)\par  The subfamily $C_1$ is isomorphic to $CI(2,2)$, an
open set in
$\Grass(2,R_2)\cong \Grass(2,6)$, so it is irreducible.  The family $C_2$ has
fiber an open set in
$\Grass(1, R_3/R_1I_2)\cong \Grass(1,5)$, over the variety ${\mathbb
P}_2^\vee\times
\mathbb{P}_2$ parametrizing a line union a point in ${\mathbb
P}^2$, so $C_2$ is irreducible, also of dimension eight. \par
We now show that the intersection $\overline{C_1}\cap\overline{C_2}=C_3$. Let
$A=R/I$ be in the intersection $\overline{C_1}\cap\overline{C_2}$. Then by
Lemma \ref{RZlem} the graded Betti numbers of $A$ are at least the supremum
$\beta (3)
$ of
$\beta(1),\beta(2)$ in Tables \ref{CI1table} and \ref{Cubictable}, so by Lemma
\ref{bettistratumlemma} they are $\beta (3)$.
Thus,
$I_2=\langle f_1,f_2\rangle$ has a ``linear'' relation in degree 3;
also, there is an extra degree 4 relation among the degree three
generator $f$, and
$f_1,f_2$. Thus, for  $A\in  \overline{C_1}\cap \overline{C_2}$ we have that up
to isomorphism
$I$ satisfies,
\begin{equation}\label{eqintideal}
I= (x V,f=hv, m^4),  V\subset R_1, \dim V=2, v\in V, h\in
R_2. 
\end{equation}
We may rule out $x\mid f$ since, as in the proof of Lemma
\ref{bettistratumlemma}, this would imply that $A$ is not level. Likewise,
if
$x\in V$  then
$v\notin
\langle x\rangle$. Hence, up to isomorphism (after a change of basis for
R) we have,  
\begin{align}\label{intersecteq}
I&= (xy,xz,yh,m^4, h\in\langle y^2,yz,z^2\rangle \text { or }\notag\\
I&=(x^2,xy,yh,m^4, h\in\langle y^2,xz,yz,z^2\rangle .
\end{align}
We now show that all ideals of the form \eqref{intersecteq}, thus all
in $C_3$ lie in the intersection
$\overline{C_1}\cap\overline{C_2} $. First, to deform those of form
\eqref{intersecteq} into $C_2$ is easy, as we need only deform
$yh$ to a degree three form having no linear factor.  To show there is a
deformation to $C_1$, we note that for $I\supset (xy,xz)$ and general
enough satisfying the first case of \eqref{intersecteq} we may assume
after a further change of variables in R (see Lemma
\ref{bettistratumlemma}), that
\begin{equation}\label{stdformeq}
I=(xy,xz,yz(y+z), m^4).
\end{equation}
We now consider the one-parameter family of ideals in $R[t]_{t}$,
\begin{equation}
I(t)= (xy,xz+tz(y+z),m^4), t\not=0.
\end{equation}
For each $t\not= 0$, we have $I(t)\in C_1$, as $(xy,zx+zt(y+z))$ is a CI.
Note that each $I(t), t\not= 0,$ contains
\begin{equation*}
 \frac{1}{t}\left[ y(xz+tz(y+z))-z(xy)\right] =yz(y+z),
\end{equation*}
hence the limit
\begin{equation}
I(0): =\lim_{t\rightarrow 0}I(t)=(xy,xz,yz(y+z),m^4).
\end{equation}
Thus, we have shown that the general element \eqref{stdformeq} of
$C_3$ is in $\overline{C_1}$ and in $\overline{C_2}$. Since $C_3$ is irreducible by  Lemma \ref{bettistratumlemma}, we
have shown
$\overline{C_1}\cap\overline{C_2}=C_3$.
\par 
Since $C_1,C_2$
are each irreducible Betti strata, since $C_3$ is
also irreducible, and since
by Lemma~\ref{bettistratumlemma},
$\beta (1),\beta (2),\beta (3)$ are the only Betti sequences that occur for level
algebras of Hilbert function $H_1$, it follows that there can be no further
irreducible components of $\Levalg (H_1)$. A check of tangent space dimensions by
{\sc Macaulay} showed that the tangent space dimension is eight for general
points of
$C_1$ or
$C_2$ but nine for the algebra $A(3)$ in
$\overline{C_1}\cap\overline{ C_2}$, defined by a monomial ideal, given by
\eqref{eqintideal2}.\par We
now consider
$\Points (T_1), T_1= (1,4,8,12,12,\ldots )$, for which
$\Delta T_1=H_1$, the sequence above. Recall that
$C_1',C_2'$ lie over over $C_1,C_2\subset \Levalg(H_1)$, and that there are
monomial ideals in each of $C_1,C_2$, as well as in $\overline{C_1}\cap
\overline{C_2}$ (equations
\eqref{ciex}, \eqref{lineeq},\eqref{eqintideal2}); thus each of
$C_1',C_2',\overline{C_1'}\cap \overline{C_2'}$ is nonempty, and the graded Betti numbers
for generic elements of these subsets agree with those that occur for the
corresponding subfamilies of $\Levalg(H_1)$ (see Tables
\ref{CI1table},
\ref{Cubictable}, and \ref{Intersecttable}). That $C_1',C_2'$ arise as stated in the Theorem is evident.
\par We first show that each of the two irreducible subfamilies
$C_1',C_2'$ of
$\Points(T_1)$ has dimension 28. Let $S=k[x,y,z,w]$, and suppose that $J\subset
S$ is the defining ideal of the punctual scheme. The dimension count for
$C_1'$ is
$\dim C_1'=16+12=28$ as $I_2=\langle f_1,f_2\rangle\in \Grass (2,10)$, has
dimension
$16=2\cdot 8$, and  12 points are chosen on the curve $\{f_1=0\}\cap
\{f_2=0\}$.
We have $\dim C_2'=28$ as follows: 
$3$ for the choice of a plane, and $4=\dim\Grass(2,4)$ for the choice of a
line in
$\mathbb P^3$, the plane union a line determining $I_2\cong \langle xy,xz\rangle$.
Then 12 for the choice of a cubic $f$ in
$S_3/(R_1I_2)$, a vector space of dimension 13 because of the linear relation
among the two generators of $I_2$. The cubic intersect the line determines three
of the 12 points; the last 9 must lie on the cubic curve defined by the cubic
surface $\{f=0\}$ intersection the plane $\{x=0\}$.\par Since a punctual
scheme must be defined by an algebra
$S/J$ having a linear nonzero divisor, it follows from our classification of
all the Betti sequences for level algebras of Hilbert function $H$
that there can be no other irreducible components of $\Points(T_1)$. 
\par
We now show that 
$C_3'=\Points_{\beta (3)}(T_1)$ has two irreducible components, the closures of $D_a'$
and $D_b'$, each of dimension 27. Evidently, $D_a'$ and $D_b'$ are
irreducible.  For $D_a'$ (9 points on the plane cubic, three on the
line), the choice of a plane cubic and a line meeting the cubic requires a
$15=3+9+3$ dimensional family, and the choice of the 12 points distributed
(9,3) gives $\dim D_a'=27$. For $D_b'$, (8 points on a plane cubic, 4 on a line), the eight
points may be chosen as a general set of 8 points on the plane, having $3+16=19$
parameters; choosing a line requires 4 parameters, and, there is one plane cubic
through the eight points meeting the line. One now chooses four points on the line
so  $\dim D_b'=19+4+4=27$. The family $D_{ab}'$ restricts one point to be the
intersection of the line and the cubic, so has dimension 26.
\par
It remains to show that there are no other components of $C_3'$. The two generators of degree two have a common factor,
 so any punctual scheme with this resolution has to lie on a plane union a line. The
twelve points can be distributed  differently between the two components, but if
there are less than eight points in the plane, the Hilbert function is at most
$7+4=11$ in degree three.  If there are at most two points on the line, the first
part of the ideal is the ideal of these points union the plane, in which case we get
socle in degree one or two. Thus, the only distributions
 that give the correct Betti numbers are the ones described in $D_a'$ and $D_b'$. 
\par
We now need to show that each point $p_\Z \in C_3'$, $p_\Z$ a punctual scheme,
lies in the closure of
$C_1'$. But, since the first part of the minimal resolution of the coordinate ring $\mathcal O_\Z$ is that of an ACM curve, one can
deform that to a complete intersection. We now show that $\overline{C_2'}\cap
\Points_{\beta (3)}(T)=\overline{D_{a}'}$. But, considering a generic point of
$D_a'$, it suffices to move the plane cubic and its nine points so that it does
not intersect the line with its three points, in order to deform to $C_2'$. On the
other hand a subfamily of ${C_2'}$ cannot converge to twelve points, of which eight are on a cubic, and four are
on a line, unless one of the points is the line intersection the cubic (so a scheme in
$D'_{ab}$).
\end{proof}
\par\noindent
{\bf Note}. Tangent space calculations using {\sc Macaulay2} \cite{GrSt} show that the
tangent space at a general point of $C_1'$,  $C_2'$ or $D_b'$ has
dimension 28, whereas the tangent space at a general point of $D_a'$ or of
$D_{ab}'$ has dimension 29.\par The strong Lefschetz (SL) condition is 
that $\times \ell$ for $\ell $ generic in $R_1$ acting on $A$ has Jordan partition $P(H)$ given by the lengths of the rows of the 
bar graph
of $H$. (This definition generalizes that in use for $H$ unimodal symmetric \cite{HW}.) It conceivably might be
 used to distinguish components of $\Levalg(H)$, as SL is an open condition. However, the algebra
$A(3)$ of
\eqref{eqintideal2} satisfies the strong Lefschetz condition, by calculation.
Since $A(3)\in C_3=\overline{C_1}\cap \overline{C_2}$, the general
elements of
$C_1,C_2$ are SL.\par
\subsection{The family $\Levalg(H_2), H_2=(1,3,6,8,9,3)$.}\label{6893}
We let $k$ be an infinite field, and consider $\Levalg (H_2), H_2=(1,3,6,8,9,3)$. The subfamily $C_1=\Levalg_{\beta (1)}(H_2)$ parametrizes type $3$, socle degree five level quotients $A$ of complete intersections
$B=R/J)$ where $J=(f_1,f_2)$ has generator degrees $(3,3)$, and having the graded Betti numbers
$\beta (1)$ given by omitting the degree four generator-relation pair from
$\beta (2)$ in Table \ref{special893table}. Thus the algebras $A=R/I$ in $C_1$
satisfy
\begin{equation}\label{comp1eq}
I_3=\langle f_1,f_2\rangle,
\end{equation}
where $f_1,f_2$ define a complete intersection in $\mathbb P^2$.
The family
$\CI(3,3)\subset \Grass(2,R_3)$ of complete intersections $B=R/(f_1,f_2)$ is
open dense in $\Grass(2,R_3)$ and has dimension
$2
\cdot 8=16$. The subfamily
$C_1$ is fibred over $\CI(3,3)$  by an open set in
the Grassmannian
$\Grass(3,R_5/R_2I_3)\cong
\Grass(3,9)$ parametrizing the choice of a type three, socle degree
five quotient of
$B=R/(I_3)$: we may take
\begin{equation}\label{comp1beq}
A=B/(W), W=\langle h_1,\ldots ,h_6\rangle \subset B_5,
\end{equation} where $W$ is a six-dimensional, general enough
subspace of the nine-dimensional space $B_5$. Thus we have
\begin{equation*}
\dim C_1=16+18=34.
\end{equation*}
One needs of course to check that such level quotients having the Hilbert
 function
$H_2$ exist: this is easy to do, using the special case $B(0)=R/(x^3,y^3)$. \par
The
graded Betti numbers $\beta(1)$ for Artinian algebras $A\in C_1$ are the minimum consistent with
the Hilbert function
$H$: there are besides the eight generators of $I$ already mentioned, ten
relations in degree six, and the three relations among the relations in degree
eight. \par The second component $\overline{C_2}$ is the closure of the Betti
stratum $C_2=\Levalg_{\beta (2)}(H_2)$ given by Table \ref{special893table}, and parametrizes ideals
whose initial portion $(I_{\le 4})$ defines a quadric union a point. More
precisely, $C_2$ consists generically of quotients
$A=R/I$ such that
\begin{align}\label{comp2aeq}
I_3&=\xi \cdot V, \xi\in R_2, V\subset R_1, \dim V =2\notag\\
I_4&=\langle R_1I_3,f\rangle 
\end{align}
and having no degree five relations. Thus $A$ is the quotient of an Artinian
algebra $B=R/(\xi\cdot V,f)$ where $R/(\xi\cdot V)$ defines a quadric
$\xi=0$ union the point $V=0$ in 
$\mathbb P^2$. 
\begin{equation}
 A=B/(W), W=\langle h_1,\ldots ,h_6\rangle \subset B_5,
\end{equation}
where as before $W$ is a six-dimensional, general enough
subspace of the nine-dimensional space $B_5$.\par
  We have 
\begin{equation*}
\dim C_2=2+5+9+18=34,
\end{equation*} 
as
follows. First, the choice of
$\xi
\cdot V, V\subset R_2$ requires 7=2+5 parameters, since the quadric requires
five and the point two parameters. 
Since
$\dim_k R_4/R_1I_3=10$, the choice of
$f\in R_4/R_1I_3$ is from an open set in $\mathbb P^9$; and the choice of
$W\subset B_5$ is that of a point in $\Grass(6,9)$.\par
 The graded Betti numbers $\beta (2)$ for an element $A\in C_2$  are those of
Table
\ref{special893table}: these graded Betti numbers are attained for
$B=R/(x^3,x^2y,z^4)$), and $W$ general enough.
\begin{table}[htb]
\begin{center}
\begin{tabular}{|l| r r r r |}
\hline
total&       1& 9 & 11  & 3 \\ \hline 
    0: &     1&   - & -  &  - \\  
     1: &  -&-&-&-\\
    2: &    - &   2& 1&    - \\
    3: &    - &   1  &  - &   -\\
    4:&     - &   6& 10&    -  \\
    5:&     - &   -& - & 3\\\hline

\end{tabular}
\caption{Graded Betti numbers $\beta (2)$ for $H_2=(1,3,6,8,9,3), \ (I_{\le
4})=(x^3,x^2y,z^4)$.}\label{special893table}
\end{center}
\end{table}
\par
There are two more Betti sequences that occur for level algebras
of Hilbert function
$H_2=(1,3,6,8,9,3)$. The first, $\beta (3)$ corresponds to ideals $I$ having an
extra relation and generator in degree five compared to the graded Betti numbers
for algebras in $C_1$. There must then be a total of seven degree five
generators. Here the ideal $J=(I_{\le 4})$ determines an algebra
$B=R/J$ defining a line union four points,
$J=(J_3)$, and
\begin{equation*}
I_3=J_3=\langle \ell g, \ell h\rangle , \, \ell\in R_1, \,\, g,h\in R_2.
\end{equation*} 
It is easy to see that this stratum has dimension 31, as follows. We count 2 for
the choice of a line
$\{\ell=0\}$, 8 for the choice of a CI $(g,h) $ of degree two forms (four points), and
(7)(3) for the choice of a three dimensional quotient $A_5$ of $R_5/R_2I_3$, a
ten-dimensional vector space.  It is also easy to see that such ideals are in the
first component $\overline{C_1}$:  writing  $J(t)=(\ell g+tf_1, \ell h+tf_2)$ we find that
$gf_2-hf_1\in J(t)$, hence is also in the limit $J(0)$, but this is no restriction, as
it merely requires $I(0)_5\cap (g,h)_5\not=0$, and a  codimension 3 space
$I_5$ always intersects a codimension four space $(g,h)_5=R_3(g,h)$ nontrivially
in the dimension 21 vector space $R_5$.\par
The last Betti stratum $\beta(4)$ is like the previous one, except for having the
extra relation and generator pair in degree four, as for algebras in $C_2$. This
corresponds to $(I_{\le 4})= (\ell g,\ell h,\ell q_3)$, where
$(g,h)$ have a common factor, a degeneration of the previous case. Thus,
\begin{equation*}
B=R/J, J=(I_{\le 4})=(\ell \rho x, \ell \rho y,\ell q_3), \ell,\rho\in R_1\quad
q_3\in (x,y)\cap R_3. 
\end{equation*}
Here it is not hard to check that this Betti stratum forms a thirty dimensional
family in the closure of the previous stratum, and as well in $\overline{C_1}\cap
\overline{C_2}$. \par
The proof below
that there are two irreducible components of
$\Levalg(H_2)$ depends primarily on the dimension count, since the poset
$\beta_{lev} (H)$ has a minimum element.
\begin{theorem}\label{36893thmA} The family $\Levalg(H_2), H_2= (1,3,6,8,9,3)$ has two
irreducible components $\overline{C_1},\overline{C_2}$, whose open dense subsets $C_1,C_2$
corresponding to two different Betti strata for $H_2$ are described above (see
\eqref{comp1beq},\eqref{comp2aeq}), each of dimension 34. Also
${\overline{C_1}}\cap C_2={V} _2$, a codimension one subset of
$C_2$ parametrizing quotients $A=R/I$,
\begin{equation}\label{specialideq}
V_2:  \quad I=( hW_1,f_4) , W_1\subset R_1, \dim_k W_1=2, h\subset R_2, f_4\in
R_3W_1.
\end{equation}
The third Betti stratum $\Levalg_{\beta (3)}(H_2)$ lies in $\overline{C_1}$ and has dimension
31. The fourth, most special Betti stratum $\Levalg_{\beta (4)}(H_2)$ is a
subscheme of
$\Levalg(H_2)$ having dimension 30, and lies in $\overline{C_1}\cap \overline{C_2}$.
\end{theorem}

\begin{proof} By Lemma \ref{RZlem} it is not
possible to specialize from a subfamily of
 $\overline{C_2}$ to a point $A$ of $C_1$, so ${\overline{C_2}}\cap
C_1=\emptyset$: this can also be seen readily since for
$B=R/J\in \overline{C_2}$,  the two-dimensional vector space $J_3$ has a common quadratic
factor, while
$I_3$ does not have a common factor for points $A=R/I$ of
$C_1$. Since the open dense $C_1\subset
\overline{C_1}$ comprises the Artinian algebras having the unique minimum possible
graded Betti numbers compatible with $H$, there can be no larger
irreducible subfamily of $\Levalg (H_2)$ specializing to both $C_1$ and
$C_2$. The family $C_2$ we have described is
exactly the subfamily having the graded Betti numbers of Table~\ref{special893table}.
Thus, that $\overline{C_2}$ has the same dimension as $\overline{C_1}$ shows that both are irreducible
components of $\Levalg(H_2)$.\par
To identify the intersection $\overline{C_1}\cap C_2$ as $V_2$ from \eqref{specialideq}, for
the moment take $W_1=\langle x,y\rangle$ and consider the following flat family
converging to a point
$B(0)$ of
$\GrAlg(H'), H'=(1,3,6,8,9,9,\ldots )$ where $B(0)=R/J(0)$ satisfies $I(0)_3$ has
a degree-two common factor $h$:
\begin{equation*}
B(t)=R/J(t),J(t)= (xh+tf_1,yh+tf_2).
\end{equation*}
 Each element of the family for $t\not=0$ contains $xf_2-yf_1$,
thus the limit $J(0)=\lim_{t\to 0} J(t)$ satisfies $ J(0)=(xh,yh, xf_2-yf_1)$.
Here $f_1,f_2$ may be chosen arbitrarily. The corresponding subfamily of $\overline{C_1}$
is fibred by an open dense subset in $\Grass (6,R_5/J(t)_5)$ over $B(t)$, parametrizing
a vector space $V(t)$, such that $I(t)=(J(t),V(t))$ defines $A(t)=R/I(t)\in \overline{C_1}$.
This shows that each element of the right side of \eqref{specialideq} occurs in
the closure of $C_1$, and vice versa.  That $V_2$ has codimension one in $C_2$
is a consequence of $``{f_4\in R_3W_1}"$ being a codimension one condition, as
$R_3W_1$ has vector space codimension one in $R_4$.
\end{proof}\par\noindent
{\bf Note.} The general element of $\Levalg_{\beta(4)}(H_2)$ is strong Lefschetz by calculation, and it follows here
as for $H_1$ that the general elements of $C_1,C_2$ are also strong Lefschetz.\par\medskip
We now consider the sum function $T_2=(1,4,10,18,27,30,30,\ldots )$ of $H_2$, and
we determine components $\overline{C_1'},\overline{C_2'}$ of $\Points (T_2)$, lying
over the components
$\overline{C_1},\overline{C_2}$, respectively, of $\Levalg (H_2)$.  Here an open dense subset $C_1'$ of
$\overline{C_1'}$ is comprised of thirty distinct points in $\mathbb P^3$ that
are general enough, lying on a CI curve $C$ that is the intersection of two
cubic surfaces.  An open dense subset $C_2'$ of
$\overline{C_2'}$ is constructed as follows:  intersect the union of a quadratic
surface
$K$ and a line $L$ with a general enough quartic surface $Q$. Choose twenty six
general enough distinct points on the curve $K\cap Q$, in addition to the four
points comprising
$L\cap Q$. The proof of the following Theorem is simplified by the fact that 
 $ H_2$ is the $h$-vector of 30 points in generic position lying on the intersection
of two cubics. \par
We denote by $S=k[x,y,z,w]$, the coordinate ring of $\mathbb P^3$.
\begin{theorem}\label{36893thmB}
The scheme $\Points(T_2), \Delta T_2=H_2=(1,3,6,8,9,3)$ has two irreducible components,
$\overline{C_1'}$, $\overline{C_2'}$, each of dimension 66. Their tangent spaces
each have dimension 66.
\end{theorem}
\begin{proof} The putative components have open dense subsets that are Betti
strata, and $\beta (1)< \beta (2)$ in $\beta_{lev} (H_2)$, so it suffices to verify the
dimension calculations. Here $CI(3,3)$ has dimension 2(18)=36, and the choice of
thirty points on the CI curve gives a total of 66 for $\dim C_1'$.\par For
$\overline{C_2'}$, the choice of a line in $\mathbb P^3$ is that of a point in
$\Grass(2,S_1)$ gives 4 dimensions, and the choice of a quadric surface
$C: c=0, c\in S_2$ up to $k^\ast$ multiple gives 9 more. The choice of a quartic
surface
$Q: q=0$ in the twenty eight dimensional vector space
$[S/(\ell_1 c,\ell_2 c)]_4,
\ell_i\in S_1$ gives 27; the choice of 26 further points on $Q\cap C$ gives a
total of 66. This shows that $\overline{C_1}$ and $\overline{C_2}$ are components.\par
In both cases, calculations in {\sc Macaulay2} 
 show that the tangent space at a general point of $C_1'$ or of $C_2'$ has dimension 66. 
\end{proof}

\subsection{An infinite sequence of examples of type three}\label{infsec}

We will now show that the example $H_2=(1,3,6,8,9,3)$ is the first in an 
infinite sequence of examples where also the number of components gets 
arbitrarily large. The idea is to start with the Hilbert function we get by 
taking a general level quotient of type three and socle degree $2c-1$ of a
 complete intersection of type $(c,c)$, where $c\ge 3$. We will see that we can 
get other level algebras with the same Hilbert function, which are not
 specializations of the ones coming from complete intersections. In order 
for us to use a result on the uniform position property \cite{Harris}, 
we have to assume in this section that $k$ is algebraically closed of characteristic $0$. 

\begin{definition}
For $c\ge 3$, we define the Hilbert function $H(c)$ of socle degree $j=2c-1$ by 
\begin{equation}\label{Hceqn}
H(c)_i = \min\{ r_i-2r_{i-c},3r_{2c-1-i}\}, \qquad 0\le i \le 2c-1,
\end{equation}
where $r_i$ is the Hilbert function of the polynomial ring $R=k[x,y,z]$. 
The transition between the two phases occurs in degree $\alpha(c):=2c-\sqrt{(c^2-1)/2}$. Note that this is usually not an integer (cf. Remark~\ref{endrmk}).
\end{definition}

We now look at the different possibilities for the generators in degree $c$ of an ideal with Hilbert function $H(c)$. 

\begin{definition}
For $c\ge 1$, let $G(c)=\Grass(2,R_c)$ parametrize two-dimensional subspaces of the space $R_c$ of forms of degree $c$ in $R=k[x,y,z]$. 
For $a=0,1,2,\dots,c-1$ let $G(c)_a$ be the subset parametrizing subspaces $V\subseteq R_c$ where the greatest common divisor has degree $a$. 
This is a stratification of $G(c)$ into disjoint semi-closed subsets.
$$
G(c) = \bigcup_{a=0}^{c-1} G(c)_a.
$$
\end{definition}

The open stratum $G(c)_0$ corresponds to vector spaces $V$ which define a complete intersection of type $(c,c)$ in $\mathbb P^2$. 
The following lemma tells us about the complement of this open stratum. Recall that $r_a=\dim_k R_a$.

\begin{lemma}\label{StrataGrass}
The dimension of the stratum $G(c)_a$ is  $r_a+2r_{c-a}-5$, for $0\le a<c$ and the closures of the smaller strata, 
$\overline{G(c)_1},\overline{G(c)_2},\dots,\overline{G(c)_{c-1}}$ are the irreducible components of the variety $X(c)$
defined as the complement of the open stratum, $G(c)_0$, in $G(c)=\Grass(2,R_{c})$. 
\end{lemma}

\begin{proof}
In the closure of the stratum $G(c)_a$ there is always a common factor of degree $a$, but there could also be a common factor
 of higher degree. In that case, this common factor has a factor of degree $a$. We get a surjective map 
$$
\mathbb P(R_a)\times \Grass(2,R_{c-a})\longrightarrow \overline{G(c)_a}
$$
by sending $(f,\langle g,h \rangle)$ to $\langle fg,fh\rangle$ and we get the dimension of $ G(c)_a$ 
as $G(c)_a= r_a-1+2(r_{c-a}-2)=r_a+2r_{c-a}-5$.\par
The variety $X(c)\subset \Grass(2,R_c)$, defined as the complement of the open stratum $G(c)_0$, parametrizes two-dimensional
 subspaces with a common factor. If the greatest common divisor of such a space $V\subseteq R_c$ is irreducible of degree $a$, 
we are in the stratum $G(c)_a$ but not in the closure of any of the other strata $G(c)_i$, $i\ne 0,a$. Thus the closures of the 
smaller strata are not contained in each other, but their union is all of $X(c)$. 
\end{proof}

Now, we will show that the we can get level algebras with Hilbert function $H(c)$ as quotients of ideals generated by two forms in degree $c$. 

\begin{lemma}\label{lemmaHF}
Let $\langle f,g\rangle\subseteq R_c$ be a general element in $G(c)_a$ and, if $a>0$, let $I=(f,g,h)$, where $h$ is a general form of degree $2c-a$.

If $a=0$ or $\sqrt{(c^2-1)/2}\le a<c$, then a sufficiently general type three level quotient of $R/I$ with socle in degree $2c-1$ has Hilbert function $H(c)$. 
\end{lemma}

\begin{proof}
We first look at the quotient of $R/I$ with the ideal $(x,y,z)^{2c}$. With $\overline I=I+(x,y,z)^{2c}$ we get that $R/\overline I$ is a level algebra with Hilbert function 
$$
H(R/\overline I)_i = \left\{\begin{array}{rl} r_i-2r_{i-c},&0\le i<2c,\\ 0,&\hbox{otherwise.}\\\end{array}\right.
$$
In particular, we have that $H(R/\overline I)_{2c-2}=H(R/\overline I)_{2c-1}=c^2$. The ideal $I$ is clearly uniquely determined by 
$\overline I$. In the parameter space of algebras having the same Hilbert function as $\overline I$, we can specialize to an ideal 
$\overline{I'}$ and if the general level quotient of $R/\overline{I'}$ has the desired Hilbert function, so does the general level quotient of $R/\overline I$.

Since $\langle f,g\rangle$ is assumed to be general, the zero-set defined by the ideal $(f,g)$ consists of a curve of degree $a$ 
and a complete intersection of type $(c-a,c-a)$. 
We can specialize $I$ into $I'=(f,g,h')$ by asserting that $h'$ vanishes on the complete intersection. This is a condition of 
codimension $(c-a)^2$. Now, in order to show that the Hilbert function of a general type three level quotient of $R/I$ is $H(c)$, 
it suffices to show that this is the case for the specialized ideal $I'$, which is an ideal of a set of $c^2$ points in the plane.
 This set of points is the union of a complete intersection of type $(a,2c-a)$ and a complete intersection of type $(c-a,c-a)$. 

We will now go on to show that if these complete intersections are chosen general enough, the Hilbert function of a general level
 quotient of type $3$ and socle degree $2c-1$ is the expected $H(c)$. 

Let $N=\max \{H(c)_i|0\le i\le c\}$ and consider a set of $N$ points in the plane which is a subset of the union of a complete 
intersection of type $(a,2c-a)$ and a complete intersection of type $(c-a,c-a)$.  We want to partition this set into three parts
 of sizes differing by at most one, such that each of the parts has the Hilbert function of a generic set of points. 

Let $W$ be a complete intersection in $\mathbb P^2$ of type $(c-a,c-a)$ having the uniform position property, i.e., all its 
subsets of the same cardinality have the same Hilbert function. Let $W=W_1\cup W_2\cup W_3$ be a partition such that $|W_1|\le |W_2|\le |W_3|\le |W_1|+1$.

For any positive integer $m< r_a$ let $\mathcal X_m\subseteq \Hilb^{m}(\mathbb P^2)$ be the subscheme parametrizing reduced sets
 of $m$ points with the generic Hilbert function, i.e., $h_i=\min\{m,r_i\}$, for $i\ge 0$.
 Let $\mathcal Y$ be the subscheme of $\Hilb^{a(2c-a)}(\mathbb P^2)$ parametrizing complete intersections of type $(a,2c-a)$ having 
the uniform position property. This is an open dense subset in the parameter space of complete intersections \cite{Harris}. 

Consider the correspondence given by 
$$
\mathcal C_m = \{(X,Y)|X\subseteq Y\} \subseteq \mathcal X_m \times \mathcal Y
$$
with the two maps ${\mathcal C}_m\longrightarrow \mathcal X_m$ and ${\mathcal C}_m\longrightarrow \mathcal Y$ induced by the projections.

The fibers of the map ${\mathcal C}_m\longrightarrow\mathcal X_m$ are all irreducible of the same dimension since we only have to pick
 general forms of degree $a$ and $2c-a$ in the ideal of $X$. Hence ${\mathcal C}_m$ is irreducible. The map $\mathcal C_m\longrightarrow \mathcal Y$ 
is finite, since a reduced complete intersection only has a finite number of subschemes. Let $\mathcal Z_m\subseteq \mathcal X_m$
 be the subscheme parameterizing sets of points, $X\subseteq \mathbb P^2$ such that the Hilbert function of one of the sets
 $X\cup W_1$, $X\cup W_2$ and $X\cup W_3$ is not generic. Then $\mathcal Z_m\subseteq\mathcal X_m$ is a closed proper subset,
 since a generic set of points in the plane yields a generic Hilbert function for each of the three sets  $X\cup W_1$, $X\cup W_2$ and $X\cup W_3$.
 Thus the dimension of the inverse image of $\mathcal Z_m$ in ${\mathcal C}_m$ is less than the dimension of ${\mathcal C}_m$ and 
hence the general element $Y$ in $\mathcal Y$ has the property that all of its subschemes of length $m$ lie outside of $\mathcal Z_m$. 
 Since this is true for all integers $m$, we can now take such a general element $Y$ in $\mathcal Y$ and let $V\subseteq Y$ be a 
subset of size $N-(c-a)^2$. We can partition $V$ into subsets $V_1\cup V_2\cup V_3$, where $|V_1|\ge |V_2|\ge |V_3|\ge |V_1|-1$. 
In this way we have a subscheme $U=V\cup W\subseteq \mathbb P^2$ such that $U_1=V_1\cup W_1$,  $U_2=V_2\cup W_2$ and $U_3=V_3\cup W_3$
 all have generic Hilbert functions. In order for this to work, we have to make sure that the sizes of the sets $V_1$, $V_2$ and $V_3$ are all less than $r_a$. It suffices to show that $N\le 3r_a-3$. 

We can get an upper bound for the number of points $N$ by looking at the value of the Hilbert function $H(c)$ in degree $\alpha(c)$ where we get
$$
N\le H(c)_{\alpha(c)} =r_{\alpha(c)}-2r_{\alpha(c)-c}= \frac{3(c^2-1)+3\sqrt{2(c^2-1)}}{4} < \frac{3}{4}(c+1)^2,
$$
for $c\ge 1$. By assumption, $\sqrt{(c^2-1)/2}\le a<c$, so we have that 
$$\begin{array}{rl}
r_a &\displaystyle =  
\frac{a^2+3a+2}{2} > \frac{c^2-1}{4} + \frac{3(c-1)}{2\sqrt2} +1 \ge \frac{c^2-1}{4}+(c-1)+1  \\&\displaystyle=
\frac{(c+1)^2+2(c-1)}{4}.\\
  \end{array}
$$
$$
r_a =
\frac{a^2+3a+2}{2} > \frac{c^2-1}{4} + \frac{3(c-1)}{2\sqrt2} +1 \ge \frac{c^2-1}{4}+(c-1)+1 =\frac{(c+1)^2+2(c-1)}{4}.
$$
Comparing these two inequalities, we get $N \le 3r_a-3$, whenever $c\ge 3$. 

The Hilbert function of the set $U$ is given by $r_i-2r_{i-c}$ in degrees $i\le \alpha(c)$ and $N$ in degrees $i>\alpha(c)$. 
\par
Suppose that an ideal $I\subset R$ satisfies $I_j\not= R_j$. We may construct a ``general'' Gorenstein ideal $J\supset I$ with $R/J$ Artinian of socle degree
$j$ by choosing first a general enough codimension one vector subspace $J_j\subset R_j$ satisfying $J_j\supset I_j$; and then taking $J$ to be the largest ideal satisfying
\begin{equation*}
J\cap M^j=(J_j)\cap M^j.
\end{equation*}
 Thus, $J$ is the \emph{ancestor ideal} of $J_j$ \cite{I3,IK}.  Equivalently, we choose a generic element $w\in (I_j)^\perp\cap \mathcal R_j$ -- that is, $w$ is annihliated
by the contraction action of $I_j$ on $R_j$ --
and let $J=\{f\in R\mid f\circ w=0\}$.\par
Now let $J(1)$, $J(2)$ and $J(3)$ be general Gorenstein ideals containing $I(U_1)$, $I(U_2)$ and $I(U_3)$, respectively, and whose quotients $R/J(1),R/J(2),R/J(3)$ each
are Artinian of
socle degree $2c-1$. 
By \cite[Lemma 1.17 or Theorem 4.1A]{IK} we have that the Hilbert function of $R/J(\ell )$ is given by 
\begin{equation*}
\max \{r_i, |U_\ell|,r_{2c-1-i}\}.
\end{equation*}
 In degrees where the Hilbert function of $U$ equals $N$, the coordinate ring of $U$ is a direct sum of the coordinate rings of $U_1$, 
$U_2$ and $U_3$. Thus the intersection $J=J_1\cap J_2\cap J_3$  gives a level algebra $R/J$ whose Hilbert function, in these degrees, 
is the sum of the Hilbert functions of the Gorenstein quotients $R/J_1$, $R/J_2$ and $R/J_3$. The sum is equal to $3r_{2c-1-i}$ 
as long as this number is less than or equal $N$. In degrees at most $\alpha(c)$, the Hilbert function of $R/J$ will equal the 
Hilbert function of $U$, i.e., $r_i-2r_{i-c}$ since the initial degrees of the three Gorenstein ideals in the coordinate rings of 
the parts are higher than $\alpha(c)$. Thus we have shown that the Hilbert function of $R/J$ is $H(c)$ and hence the general level
 quotient has Hilbert function $H(c)$. 
\end{proof}

After having established that there are level algebras with Hilbert function $H(c)$ with different degrees of the common divisor in degree $c$, 
we can now state the main theorem of this section.

\begin{theorem}\label{InfiniteTypeThree}
For $a=0$ and $\sqrt{(c^2-1)/2}\le a<c$, the family $\mathcal F_a$ of all level
algebras in $\Levalg(H(c))$ whose degree $c$ part lies in $G(c)_a$ is a non-empty open set in a component of $\Levalg(H(c))$. 
The dimension of this component is $4c^2+3c-11$ for $a=0$ and $4c^2 + 3c-12+(c-a)^2$ otherwise. $\Levalg (H(c))$ is reducible for $c\ge 3$.
\end{theorem}

\begin{proof}
Consider the map 
$$
\Phi:\Levalg(H(c))\longrightarrow \Grass(2,R_c)
$$
given by sending the level ideal $I$ to its degree $c$ component $I_c$. 

By Lemma~\ref{lemmaHF} we know that $\mathcal F_a$ is non-empty, since $\Phi$ has a non-empty fiber over the general point of $G(c)_a$.
 The general point in $\Grass(2,R_c)$ lies in the image and we conclude that $\Phi$ is dominant. There has to be at least one component 
of $\Levalg(H(c))$ dominating the image, and since the general fiber of $\Phi$ is irreducible of dimension $3(c^2-3)$, we conclude that
 there is a single such component, $\overline{\mathcal F_0}$, of dimension
$$
2r_c-4 + 3(c^2-3) = 4c^2+3c-11.
$$
Let $X(c)$ be the variety defined as the complement of $G_0$ in $\Grass(2,R_c)$. Lemma~\ref{StrataGrass} gives us the dimensions of the components, 
$\overline G(c)_a$, of $X(c)$. Now we know from Lemma~\ref{lemmaHF} that for $\sqrt{(c^2-1)/2}\le a<c$, the fiber over the general point of $G(c)_a$ is non-empty. 
Thus we know that there has to be a component of $\Levalg(H(c))$ dominating the component $\overline G(c)_a$ of $X(c)$. The fiber of $\Phi$ over the general point in $G(c)_a$ is 
irreducible of dimension $r_{2c-a}-2r_{c-a}+3c^2-9$ since we have to add a form of degree $2c-a$ to the ideal generated by the two forms
 of degree $c$ due to their common factor of degree $a$. Observe that since $a\ge \sqrt{(c^2-1)/2}$ we have that $2c-a\le \alpha(c)$ and the form of degree $2c-a$ 
is determined up to scalar multiples by the level ideal. Thus the component, $\overline{\mathcal F_a}$, dominating $\overline G(c)_a$ is irreducible of dimension
$$
2r_{c-a}+r_a-5 + r_{2c-a}-2r_{c-a} + 3c^2-9 = r_a+r_{2c-a} + 3c^2-14 = 4c^2+3c-12+(c-a)^2.
$$
Since this is greater than or equal to the dimension of the component spanned by $\mathcal F_0$, we have that $\mathcal F_a$ cannot be 
contained in the closure of $\mathcal F_0$, and we have shown that they are different components. For different $a\ne 0$, the components 
spanned by $\mathcal F_a$ are different since they map onto different components of $X(c)$. 
\end{proof}

\begin{remark}\label{endrmk}
The first few Hilbert functions in the series $H(c)$ are 
$$
\begin{array}{cccccccccccccccc}
H(3)&=(1,&3,&6,&8,&9,&3)\\
H(4)&= (1,&3,&6,&10,&13,&15,&9,&3)\\
H(5)&=(1,&3,&6,&10,&15,&19,&22,&18,&9,&3)\\
H(6)&=(1,&3,&6,&10,&15,&21,&26,&30,&30,&18,&9,&3)\\
H(7)&=(1,&3,&6,&10,&15,&21,&28,&34,&39,&43,&30,&18,&9,&3)\\
\end{array}
$$
and the first time we get more than two components is for $c=7$. The number of components given by Theorem~\ref{InfiniteTypeThree}
 is $c+1-\lceil\sqrt{(c^2-1)/2}\rceil$, which is bounded from below by $(1-1/\sqrt2)c$. 

The values of $c$ for which the degree $\alpha(c)=2c-\sqrt{(c^2-1)/2}$ is an integer correspond to solutions of Pell's equation, 
$c^2-2d^2=1$, which in turn are related to the continued fraction expansion of
$\sqrt 2$, namely $2= (1,2,2,\ldots )$.  Considering the approximants $p_k/
q_k$, and letting $c=p_k, d=q_k$ for $k$ odd, one obtains all such solutions.
The first five are $(p_k,q_k)= (3,2),(17,12), (99,70), (577, 408),
(3383,2378),$ for $k=1,3,5,7,9$. (See \cite[Theorem 13.11]{Ro}).\par
It would be interesting to know whether the infinite series of Artinian examples of Theorem~\ref{InfiniteTypeThree}
 lifts to an infinite series of examples with level set of points in a similar way as in the first example in the series.
 This question is still open and the problem
 is that we are lacking general results that guarantee that the Hilbert function 
of a general enough level set of points on a given curve has the expected Hilbert function.

\end{remark}

\begin{ack} The authors are grateful for the comments and interest of several
in this work, in particular J. Migliore and J.-O. Kleppe. The authors are grateful
 to the Fields Institute, and organizers Sara Faridi et al of their January 2005 Conference
  ``Resolutions, Inverse Systems, and Coinvariants", held at the University of Ottawa 
for its interest and the occasion to meet; and also to the conference on Algebraic Geometry
 in August, 2005 at Seattle.  The second author is grateful for an invitation to  visit KTH 
in April, 2005, which was made possible by a grant from the G\"oran Gustafsson Foundation.  
The authors are grateful for the occasion of the MAGIC Conference at Notre Dame, of 2005, 
organized by A. Corso, J. Migliore, and C. Polini. An early version of the article 
was submitted to the Proceedings of the MAGIC conference, and we appreciate referee
 comments including the suggestion to include a proof of an announced result
 concerning an infinite series of examples. This would have changed the submission 
significantly, and we have instead written the present article.  We gladly acknowledge 
our debt to the MAGIC organizers, and to the referees.
\end{ack}

\bibliographystyle{amsalpha} 

\end{document}